\documentclass[12pt,onecolumn]{IEEEtran}

\IEEEoverridecommandlockouts                              

\usepackage{float}
\usepackage{graphics}
\usepackage{subfigure}
\usepackage{graphicx} 
\usepackage{epsfig} 
\usepackage{amsmath} 
\usepackage{amssymb}  

\usepackage{mathdots}
\usepackage{MnSymbol}
\usepackage[ruled]{algorithm2e}
\usepackage{amsmath}
\usepackage[usenames,dvipsnames,svgnames,table]{xcolor}

\newtheorem{theorem}{Theorem}
\newtheorem{remark}{Remark}
\newtheorem{lemma}{Lemma}
\newtheorem{definition}{Definition}
\newtheorem{proposition}{Proposition}
\newtheorem{corollary}{Corollary}

\newcommand\g{\cellcolor{gray!10}}
\newcommand\gd{\cellcolor{gray!20}}
\usepackage{tikz}
\newcommand\myline[1][]{%
  \,\tikz[baseline]\draw[very thick,#1](0,-\dp\strutbox)--(0,\ht\strutbox);\,%
}

\title{
Selective Strong Structural Minimum Cost Resilient Co-Design for Regular Descriptor Linear Systems
\vspace{-0.2cm}
}

\author{Nipun Popli $^{\star,\sharp}$ \ S\'ergio Pequito $^{\diamond,\sharp}$ \ Soummya Kar $^{\star}$  \  A. Pedro Aguiar $^{\dag,\ddag}$ \ Marija Ili\'c $^\star$
\thanks{ This work was partially supported by the CMU-Portugal (ICTI) program, and by projects NSF CIF-1513936, POCI-01-0145-FEDER-006933/SYSTEC funded by FEDER funds through COMPETE2020 and FCT.}
\thanks{$^\sharp$ Both authors contributed equally to this work.}
\thanks{
$^{\star}$ Department of Electrical and Computer Engineering, Carnegie Mellon University, Pittsburgh, PA 15213.}
\thanks{$^{\diamond}$ Department of Electrical and Systems Engineering, School of Engineering and Applied Science, University of Pennsylvania, Philadelphia, PA 19104, USA}
\thanks{
$^{\dag}$ Institute for System and Robotics, Instituto Superior T\'ecnico,  University of Lisbon, Lisbon, Portugal.
        }%
\thanks{
$^{\ddag}$ Department of Electrical and Computer Engineering, Faculty of Engineering, University of Porto (FEUP), Porto, Portugal  }%
\vspace{-0.9cm}
}

\begin{document}

\maketitle
\thispagestyle{empty}
\pagestyle{empty}

%
%

\begin{abstract}

This paper addresses the problem of minimum cost resilient \mbox{actuation-sensing-communication} co-design for regular descriptor systems while ensuring selective strong structural system's properties. More specifically, the problem consists of determining the minimum cost deployment of actuation and sensing technology, as well as communication between the these, such that decentralized control approaches are viable for an arbitrary realization of regular descriptor systems satisfying a pre-specified selective structure, i.e., some entries can be  zero, nonzero, or either zero/nonzero. Towards this goal, we rely on strong structural systems theory and extend it to cope with the selective structure that casts resiliency/robustness properties and uncertainty properties of system's model. Upon such framework, we introduce the notion of selective strong structural fixed modes as a characterization of the feasibility of decentralized control laws. Also, we provide necessary and sufficient conditions for this property to hold, and show how these conditions can be leveraged to determine the minimum cost resilient placement of \mbox{actuation-sensing-communication} technology ensuring feasible solutions. In particular, we study the minimum cost resilient actuation and sensing placement, upon which we construct the solution to our problem. Finally, we illustrate the applicability the main results of this paper on an electric power grid example.
\end{abstract}

\vspace{-0.3cm}

%
%

\section{Introduction}\label{intro}
Over the past decade, we have witnessed a steady growth of large-scale systems, and examples include electric power grid~\cite{ilic2000dynamics}, transportation networks~\cite{vsiljak2008dynamic}, biological~\cite{Dai89} and social networks~\cite{liu11}, and swarms of multi-agents~\cite{TannerPK04}, just to name a few. In fact, it is often necessary to evaluate the control theoretic properties of such systems, e.g.  controllability and observability, which are critical to ensuring the systems' proper dynamical evolution~\cite{vandeWal2001487}. Also, more than often, the large-scale and geographical nature of such systems entails decentralized data sharing with the actuators. Consequently, the data accessible to each actuator must be sufficient to ensure the closed-loop system specifications~\cite{Sandel_etall:1978}. The \emph{information pattern} captures the information accessibility to actuators and implicitly represents the communication requirements associated with the decentralized control scheme. It is equally important to note that while dealing with such large-scale systems, \mbox{actuation-sensing-communication} has to be simultaneously designed to ensure the existence of decentralized solutions~\cite{Skogestad2004219}. Also, due to the infrastructure and maintenance costs, it may be desirable to identify the minimum \mbox{actuation-sensing-communication} required to ensure system specifications. Besides, as a consequence of the susceptibility of the large-scale systems to component failures, their specifications often have to address robustness/resilience properties. Therefore, the objective is to ensure that the control properties hold in the case of \mbox{actuation-sensing-communication} failure, or compromised due to an external agent, while incurring minimum investment cost.

Furthermore, parametric uncertainties in the system model are inevitable, and even if that is not the case, assessment of control theoretic properties based on numerical methods is impractical when the dimension of the system is large~\cite{Reinschke:1988}. Therefore, in this paper, we aim to characterize the minimum \mbox{actuation-sensing-communication} for all possible realizations of the linear system plant satisfying a given structure that captures the interconnections between different assets of the dynamical system. Toward this goal, we rely on {strong structural systems theory}~\cite{Jarczyk11}, which aims to ensure control theoretic properties for all possible nonzero realizations of a given structure, and  we extend this to account \emph{for all} possible realizations obeying a given structure that identifies which entries are strictly zero, strictly nonzero, or possibly zero and nonzero, which we refer to as \emph{selective} strong structural systems. Notice that it is distinct from the structural systems properties which only ensure that almost all parameterizations guarantee the control theoretic properties~\cite{Dion03survey}. Nonetheless, when dealing with interconnected dynamical systems, it may occur that, in practice, structural properties do not hold, which motivates the need for strong structural systems theory, whereas the fragility of some interconnections or their small dynamical dependency prompts the need for a selective strong structural systems approach proposed in this paper. Also, we intend to use selective strong structural systems in the general context of regular descriptor linear time-invariant systems, which account for scenarios commonly found across different interconnected systems with conservation laws, for instance, in an electric power grid~\cite{LuenbergerDescriptorSystems}.
{
Contrarily to linear time-invariant systems, there have been proposed
a variety of possible definitions for controllability and observability (see~\cite{berger2013controllability}). Therefore, when referring to these concepts, we adopt the definition that is closest to the one used in linear time-invariant systems, which is commonly associated with the state reachability. For instance, by \emph{controllability} we mean \emph{R}-controllability~\cite{Yip}, which is often also referred to as behavoral controllability~\cite{berger2013controllability}.}

The interplay between structure of the system and its specific parametric descriptions are the scope of {structural}~\cite{Dion03survey} and strong structural systems theory~\cite{Jarczyk11}. Although a considerable amount of work has addressed structural systems properties and \mbox{actuation-sensing-communication} co-design (see~\cite{PequitoC18,PequitoJ1,clark2016submodularity} and references therein), the same is not true regarding strong structural systems properties. 
The notion of strong structural controllability was introduced in~\cite{Mayeda79},  and necessary and sufficient conditions for linear time-invariant systems were provided in~\cite{Reinschke92,Jarczyk11}, as well as for linear time-varying systems in~\cite{Reissig14}; in particular, the necessary and sufficient conditions for linear  time-varying systems can be evaluated in terms of an auxiliary linear time-invariant system~\cite{Reissig14}.  An interesting and pedagogical example of the applicability of strong structural systems can be found in~\cite{Bowden12}. As an alternative to the state space representation, the strong structural systems theory has also been proposed to study systems properties in the frequency domain~\cite{Daasch201651}.  Whereas the problem of verifying strong structural controllability can be addressed in linear-time complexity~\cite{Weber14}, the problem of selecting the minimum number of actuators out of a possible set of configurations was shown to be in general NP-hard~\cite{Chapman13,Trefois15}.  In~\cite{Monshizadeh14}, a graph-theoretic characterization of strong structural controllability is provided in the context of leader-follower, and later extended to account for the study of target controllability, i.e., the controllability of a subset of state variables~\cite{van2016distance}.   In~\cite{camsap}, the authors have introduced new algebraic necessary and sufficient conditions for strong structural controllability for linear time-invariant systems, and studied the  minimum placement of \emph{dedicated actuators}, i.e., actuators that manipulate a single state variable. In the present paper, we extend these results to the case where the sparsest solutions are sought in the context of regular descriptor systems and selective strong structural systems. Furthermore, we address the minimum cost resilient \mbox{actuation-sensing-communication} co-design problem, which necessitates the introduction of novel concepts in strong structural theory, as well as necessary and sufficient conditions to enable the design. Finally, we notice that the techniques used to address the co-design problem in the context of strong structural systems are algebraic, and, therefore, quite different from the graph theoretic conditions used in structural systems theory~{\cite{PequitoJ10,PequitoC18,PequitoJ1}}.


The main contributions of this technical note are threefold: (\emph{i}) we introduce the concept of selective strong structural fixed modes that ensures the non-existence of fixed modes to any realization of the system's descriptor state space representation satisfying a specified structure; (\emph{ii}) we address the minimum cost resilient \emph{co-design} of \mbox{actuation-sensing-communication} that ensures the non-existence of selective strong structural fixed modes; and in addition, (\emph{iii}) we address the sparsest actuator (respectively, sensor) design for descriptor linear time-invariant systems that ensures selective strong structural controllability (respectively, observability) upon which we build the solution to the co-design problem.

The remainder of the paper is organized as follows: Section~\ref{probl_formu} introduces the selective strong structural notion of decentralized fixed modes for regular linear time-invariant descriptor systems, and, subsequently, the formal problem statements are presented. Next, Section~\ref{main_result} begins with the review of concepts in strong structural system theory and new definitions are provided. Subsequently, based on strong structural theory, the solutions to the problem statements are presented. Section~\ref{illus_examp} illustrates the application of proposed solutions on a sixteen dimensional multi-input multi-output model of a $5$-bus electric grid. Section~\ref{concl_furth} concludes the paper and further research avenues are discussed.

%
%

\section{Problem Formulation}\label{probl_formu}
Consider a dynamical systems modeled, or locally approximated, by a { regular} descriptor system given by:
\begin{align}
E\dot x \left(t \right) &=Ax\left(t \right) + Bu\left(t \right), \label{dynLinear}\\ 
y\left(t \right) &=Cx\left(t \right),\label{outputLinear} 
\end{align}
where the state vector evolution is represented by $x\left(t\right)\in \mathbb{R}^n$, the input vector $u\left(t\right) \in \mathbb{R}^p$, and the output vector $y\left(t\right)\in \mathbb{R}^m$ over time $t \in \mathbb{R}_+$. { In addition, the dynamics by $A\in\mathbb{R}^{n\times n}$, the descriptor matrix is denoted by $E\in\mathbb{R}^{n\times n}$ such that  $\det(E-\lambda A)\not\equiv 0$  (i.e., $\det(E-\lambda A)\neq 0$ for almost all $\lambda\in\mathbb{C}$), the input matrix by $B\in \mathbb{R}^{n\times p}$, and the output matrix by $C \in \mathbb{R}^{m\times n}$.} We refer to the descriptor system in~\eqref{dynLinear}-\eqref{outputLinear} by the tuple $\left(E, A, B, C\right)$. Additionally, we can think about the system in a \mbox{closed-loop}, where one potential strategy is to use output feedback under partial information constraints. The availability of the measurements to each actuator is captured by the notion of information pattern. This can be described by a $p\times m$ binary matrix $\bar K\in\{0,\times\}^{p\times m}$, where an entry $\bar K_{i,j}=\times$ if the data from sensor $j$ is available to actuator $i$, and zero otherwise. In other words, the entry $\bar K_{i,j}=\times$ denotes the existence of a communication channel from the $j$-$\text{th}$ sensor to the $i$-$\text{th}$ actuator. Furthermore, one can consider static output feedback, where the input response is designed as a linear combination of the available measurements from the sensors, i.e.,  
\begin{equation}
u\left(t \right)=-Ky\left(t \right),\label{static_fdbk}
\end{equation}  
where $K \in \mathbb{R} ^ {p \times m}$ is the feedback gain matrix, whose sparsity is induced by the information pattern, i.e., $K_{i,j}=0$ if $\bar K_{i,j}=0$. The static output feedback can leverage the limited communication and computational capabilities in large scale dynamical systems~\cite{Syrmos}. We represent the closed-loop static output feedback descriptor system \eqref{dynLinear}-\eqref{static_fdbk}, under the information pattern constraint $\bar K$, by the tuple $\left( {E, A, B, C; \bar K} \right)$. Specifically, the communication design must ensure  the existence of feedback gain matrices $K \in \mathbb{R} ^ {p \times m}$, with the sparsity of information pattern $\bar K$, to change the static output feedback closed-loop modes. The modes that cannot be changed by such gains are known as \emph{fixed modes}~\cite{Davison, Anderson1981703}. Besides, it is well known that the controllability of the tuple $\left( {E, A, B} \right)$, and the observability of the tuple $\left( {E, A,C} \right)$, are necessary but not sufficient for existence of a control law based on a static output feedback. 


Nonetheless, the efficacy of verifying controllability, observability, or the existence of fixed modes with respect to the information pattern $\bar K$, is contingent on the numerical accuracy of the parameters in $\left( {E, A, B, C; \bar K} \right)$. Therefore, to deal with such scenarios, we propose to rely on strong structural theory~\cite{Jarczyk11}. The strong structural theory enables the investigation of basic control properties based solely on sparsity pattern of the system plant matrices in~\eqref{dynLinear}-\eqref{static_fdbk}. The matrix entries are qualitatively represented as either nonzero, denoted by $\times$, or zero. However, in a more general scenario, it may not be known whether some entries are zero or nonzero. Therefore, a framework capable of accounting for such a scenario is often desired. Hereafter, we provide such framework which we refer to as \emph{selective strong structural theory}. In particular, besides the zero and nonzero entries, we selectively allow some entries to be any real value, which we denote by $\otimes$. As a consequence, within this framework, the \emph{selective} structural matrices $\bar X \in {\left\{ {0,\times,\otimes} \right\}^{a\times b}}$ define an equivalent class of matrices as \vspace{-0.2cm}
\[\left [ \bar X \right]= \left\{ {\begin{array}{*{20}{c}}
{X{ \in {\mathbb{R}}^{a \times b}}:}&{\begin{array}{*{20}{c}}
{{\rm{if}}\;{{\bar X}_{ij}} = 0\;{\rm{then}}\;{X_{ij}} = 0,\ \ \ \ \ \ \ \ \ }\\
{\begin{array}{*{20}{l}}
{{\rm{if}}\;{{\bar X}_{ij}} =  \times \;{\rm{then}}\;{X_{ij}} \in \mathbb{R}\backslash \left\{ 0 \right\},}\\
{{\rm{if}}\;{{\bar X}_{ij}} =  \otimes \;{\rm{then}}\;{X_{ij}} \in \mathbb{R}}
\end{array}}
\end{array}}
\end{array}} \hspace{-0.2cm}\right\}.
\]
{
Also, because we consider only regular descriptor systems, we also need to introduce the following subclass for $\bar A,\bar E \in {\left\{ {0,\times,\otimes} \right\}^{n\times n}}$ :
 \[
 \left([\bar E],[\bar A]\right)^{\star}=\{(E,A)\in[\bar E]\times[\bar A] : \det(A-\lambda E)\not\equiv 0\}.
 \]
}
These equivalent classes are used to represent the descriptor systems in selective strong structural system theory, i.e., these are represented by the tuple {$\left((\left[\bar E \right], \left[\bar A \right])^{\star}, \left[\bar B \right], \left[\bar C \right], \bar K \right)$}, which, with some abuse of notation, we represent as the tuple of selective structural matrices as $\left(\bar E , \bar A , \bar B, \bar C; \bar K \right)$.

Now, we introduce some preliminary terminology required to characterize the solutions to the problems explored in this paper.  Let $M$ be a $m_1\times m_2$ matrix, then we refer to $m_1$ and $m_2$ as the \emph{height} and \emph{length} of the matrix, respectively. In addition, we need the following definitions~\cite{camsap}.

\begin{definition}[Stair matrix]
A matrix $M\in\{0,\times,$ $\otimes\}^{m_1\times m_2}$ is said to be a \emph{stair matrix} if it is of the form
\begin{equation*}
M_{m_1\times m_2}=\left[\begin{array}{cccccc}
\gd S^1_{n_h^1\times n_l^1} & \g   &\g&\g &\g \mathbf{0}_{n^1_h\times(m_2-n_l^1)} &\g\\
\gd S^2_{n_h^2\times n_l^2} &\gd & \g&\g  &\g \mathbf{0}_{n^2_1\times(m_2-n_l^2)} &\g \\
&\ddots &&&\ddots&\\
\gd S^k_{n_h^k\times n_l^k}  & \gd & \gd  &\gd &\g \mathbf{0}_{n_h^k\times(m_2-n_l^k)} &\g \\
\end{array}\right],
\end{equation*}
where each $n_h^i\times n_l^i$ matrix $S^i_{n_h^i\times n_l^i}$ denotes the $i$-th \emph{step} ($i=1,\ldots,k$) such that  $n_l^j< n_l^{j+1}$ ($j=1,\ldots,k-1$), and  $\mathbf{0}_{p_1\times p_2}$ denotes the $p_1\times p_2$ zero matrix. In addition, $M$ is in the \emph{maximal stair form}, if  there exist no permutation matrices $P^M_r$ and $P^M_c$ such that $P_r^MMP_c^M$ has more steps and zero matrices with larger length than those in $M$.\hfill $\diamond$
\label{stepmatrix}
\end{definition}

Notice that the steps in a stair matrix $M$ are ordered from top to bottom by length, i.e., from the smallest $S^1_{n_h^1\times n_l^1}$ to the largest $S^k_{n_h^k\times n_l^k}$. Now, given a stair matrix, we introduce the notion of \emph{step difference}.
 

\begin{definition}[Step difference]
Given a  stair matrix $M\in\{0,\times,\otimes\}^{m_1\times m_2}$ with $k$ steps, the $i$-th \emph{step difference} denoted by $\Delta^{i}$ can be recursively defined as follows:
\begin{enumerate}
\item[(i)] $\Delta^{1}=S^1_{n_h^1\times n_l^1}$; and
\item[(ii)] $\Delta^{i+1}=S^{i+1}_{n_h^{i+1}\times n_l^{i+1}}[:,n_l^i+1:n_l^{i+1} ]$, for $i=1,\ldots,k-1$,
\end{enumerate}
where $P[:,c_1:c_2]$ corresponds to the submatrix of a matrix $P$ comprising all the rows and columns indexed from $c_1$ to $c_2$.
\hfill $\diamond$
\label{stepdifference}
\end{definition}

Simply speaking, from Definition~\ref{stepdifference}, it follows that a step difference $\Delta^{i+1}$ is a result of the `difference' between two adjacent steps $S^i_{n_h^{i}\times n_l^{i}}$ and $S^{i+1}_{n_h^{i+1}\times n_l^{i+1}}$ ($i=1,\cdots,k-1$), in the sense that it contains the same rows of $S^{i+1}_{n_h^{i+1}\times n_l^{i+1}}$ but only the columns from $n_l^i+1$ to $n_l^{i+1}$, illustrated as follows:
 \begin{equation*}
\left[\begin{array}{cccccc}
\gd S^i_{n_h^i\times n_l^i} & \gd   &\g&\g &\g \mathbf{0}_{n^i_1\times(m_2-n_l^i)} &\g \\
\gd S^{i+1}_{n_h^{i+1}\times n_l^{i+1}}  & \gd & \hspace{-0.2cm} \myline[dashed] \  \gd\Delta^{i+1}  &\g &\g \mathbf{0}_{n^{i+1}_1\times(m_2-n_l^{i+1})} &\g \\
\end{array}\right].
\end{equation*}

{
 Let $\bar A^{\lambda}=\bar A - \lambda \bar E$, where ${\bar A^{\lambda}}_{ij}=\times$ if $\bar A_{ij}=\times$ and $\bar E_{ij}=0$, ${\bar A^{\lambda}}_{ij}=\otimes$ if either $\bar A^{\lambda}_{ij}=\otimes$ or $\bar E_{ij} \ne 0$, and $\bar A^{\lambda}_{ij}=0$ if $\bar A_{ij}=0$ and $\bar E_{ij}=0$.  In what follows, we will focus on systems that satisfy the following assumption (see Appendix for further details).

\textbf{Assumption~1:} 
Each step difference in the stair matrix of $\bar A^{\lambda}$ has one column vector such that for all  its parametric choices,  the remaining vectors in the step difference admit a parameterization that makes all vectors proportional to each other. \hfill $\circ$

}


In addition, real world systems often experience unexpected failures of actuators, sensors, communication links or their combination. Thus, the design of large-scale systems must account for such possible failures. More specifically, \mbox{actuation-sensing-communication} must be co-designed such that the different control-theoretic properties hold after occurrence of such adverse events. Despite such considerations, in real world setups, the deployment of the \mbox{actuation-sensing-communication} infrastructure incurs multitude expenditures, such as cost of devices, installments, and their maintenance. Consequently, it is often required to consider the minimum cost \mbox{actuation-sensing-communication} co-design, that guarantees certain degree of resiliency with respect to \mbox{actuation-sensing-communication} failures. Motivated by the importance of such problems in practice, together with lack of knowledge of all system's parameters, we propose to address the following three different (but related) problems.
\newpage

$\triangleright$ \underline{Minimum Cost Resilient Actuation Selection Problem}

First, we need to introduce the selective strong structural counterpart of controllability, that readily extends the notion of controllability{~\cite{Yip}} in strong structural theory as follows.

\begin{definition}\textbf{\emph{(Selective Strong Structural Controllability (SSSC))}}
The tuple $\left( {\bar E, \bar A, \bar B} \right)$ is selective strong structural controllable if and only if $\left( { E,  A,  B} \right)$ is controllable for all {$ (E,A)\in \left(\left[\bar E\right],\left[\bar A\right]\right)^{\star}$} and  $B\in \left[\bar B\right]$. 
 \hfill $\diamond$
\label{def_SSSC}
\end{definition}

In addition, different actuators can (potentially) actuate different state variables while incurring different costs. Subsequently,  a heterogenous costs or weights are associated with the actuation of the states, and can be represented using a weight matrix $W^B \in \mathbb{R}_+^{n \times p}$, where the entry $W^B_{ij}$ represents the cost of actuating the state with index $i$ by the actuator with index $j$. In other words, the actuation cost depends only on the state variable actuated and not on the actuator performing the control. Therefore, the first problem we address in this paper is stated as follows.

\noindent $\mathcal P_1$ Given selective structural matrices $\bar E, \bar A \in \{ 0, \times,$ $\otimes\}^{n \times n}$, a maximum of $k$ actuator failures, and actuation cost structure $W^B\in\mathbb{R}_+^{n \times(k+1)n}$, determine $\bar B^*$ that solves the following problem
\begin{align}
\mathop {\min }\limits_{\bar B \in {{\left\{ {0, \times } \right\}}^{n \times \left(k+1\right)n}}} & \qquad \qquad \qquad\left\| {\bar B} \right\|_{{W}^B}\\
\text{s.t.}\quad & \left( {\bar E, \bar A,\bar B\left( {{\mathcal{I}_B}} \right)} \right) \text{ is SSSC,} \text{ for all } \mathcal{I}_B \subset \mathcal{N}, \ \mathcal{I}_B' \subset \mathcal{N},  \notag \\  
& \text{with } \mathcal{I}_B = \mathcal{N} \backslash  \mathcal{I}_B',  \text{and } {\left| {\mathcal{I}_B'} \right| \le k},\notag
\end{align}
where ${\left\| {\bar B} \right\|_{{W}^B}} = {{\bf{1}}^T}\left( {\bar B \odot {W}^B} \right){\bf{1}}$ where $\bar B \odot {W^B}$ is defined as $\left[{\bar B} \odot {W^B} \right]_{ij}= W_{ij}^B$ if ${{\bar B}_{ij}}=\times$ and $0$ otherwise. In addition, $\bf{1}$ represents the ones vector with appropriate dimensions, the structural matrix $\bar B \left( {{\mathcal{I}_B}} \right)$ is the subset of columns corresponding to the actuators with indices in $\mathcal{I}_B$, $\mathcal N = \left\{ {1,2,3, \ldots ,\left(k+1\right)n} \right\}$ and $\mathcal{I}_B'$ contains the indices of the columns representing actuators that have malfunctioned. 
\hfill $\diamond$

In the minimum cost resilient actuator selection problem, notice that the entries in the selective structural matrix $\bar B$ are restricted to nonzero $\times$ and zero $0$, since any real entry $\otimes$ does not allow for posing of a well-defined cost objective in the problem. In addition, we consider an $n \times \left(k+1\right)n $ structural matrix $\bar B$ to allow a feasible solution to the problem $\mathcal{P}_1$, since the concatenation of selective structural pattern of $\left(k+1\right)$ identity matrices is granted to achieve feasibility, so that at least the solution to the problem exists. Further, we notice that $\bar B^*$ may contain \mbox{zero-columns}, and its nonzero columns will be associated with the \emph{effective actuators} that are considered in the design procedure.

$\triangleright$ \underline{Minimum Cost Resilient Sensing Selection Problem}

Similar to the previous problem, we formalize the minimum cost resilient sensing selection problem, for which we need the following definition.

\begin{definition}\textbf{\emph{(Selective Strong Structural Observability (SSSO))}} 
The tuple $\left( {\bar E, \bar A, \bar C} \right)$ is selective strong structural observable if and only if $\left( {E,  A,  C} \right)$ is observable   for all  {$ (E,A)\in \left(\left[\bar E\right],\left[\bar A\right]\right)^{\star}$} and  $C\in \left[\bar C\right]$. 
 \hfill $\diamond$
\label{def_SSSO}
\end{definition}


In addition, we define the sensing cost matrix $W^C \in \mathbb{R}_+^{m \times n}$, in which the entry $W^C_{ij}$ represents the cost of measuring the $j$-th state variable by the sensor with index~$i$. In other words, the sensing cost depends only on the state variable measured and not on the sensor performing the measurement. Subsequently, the second problem we address is posed as follows.

\noindent $\mathcal P_2$  Given selective structural matrices $\bar E, \bar A \in  \{0, \times,$ $\otimes\}^{n \times n}$,  a maximum of $k$ sensor failures, and sensing cost structure $W^C\in\mathbb{R}_+^{(k+1)n \times n}$, determine $\bar C^*$ that solves the following problem
\begin{align}
\mathop {\min }\limits_{\bar C \in {{\left\{ {0, \times } \right\}}^{\left(k+1\right)n \times n}}} & \qquad \qquad \qquad\left\| {\bar C} \right\|_{{W}^C}\\
\text{s.t.}\quad & \left( {\bar E, \bar A,\bar C\left( {{\mathcal{I}_C}} \right)} \right) \text{ is SSSO,}  \text{ for all } \mathcal{I}_C \subset \mathcal{N}, \ \mathcal{I}_C' \subset \mathcal{N}, \notag\\
& \text{with }  \mathcal{I}_C = \mathcal{N} \backslash  \mathcal{I}_C', \text{and } {\left| {\mathcal{I}_C'} \right| \le k},\notag
\end{align}
where $\bar C \left( {{\mathcal{I}_C}} \right)$ is the subset of rows corresponding to the sensors with indices in $\mathcal{I}_C$, and $\mathcal{I}_C'$ contains the indices of the rows representing sensors that have malfunctioned.~\hfill~$\diamond$

$\triangleright$ \underline{Minimum Cost Resilient Actuation-Sensing-Communication Selection Problem}

Lastly, we introduce the notion of fixed modes in the context of the selective strong structural systems theory for {regular} descriptor systems, which we refer to as \emph{selective strong structural fixed modes}, which readily extends the characterization in~\cite{Anderson1981703}.

\begin{definition}\textbf{\emph{(Selective Strong Structural Fixed Modes  (SSSFM))}}
A { regular} descriptor system with a given selective structural pattern $\left( {\bar E, \bar A, \bar B, \bar C; \bar K} \right)$ has a SSSFM  $\lambda \in \mathbb{C}$ (with respect to the information pattern $\bar K$), if there exists {$ (E,A)\in \left(\left[\bar E\right],\left[\bar A\right]\right)^{\star}$}, $B \in \left [ \bar B \right]$, $C \in \left [ \bar C \right]$ and $K$ satisfies the information pattern $\bar K$, such that 
$
rank\left( {\lambda E -A-BKC} \right)<n. 
$
 \hfill $\diamond$
\label{struct_fixedmodes}
\end{definition}

Additionally, to setup the communication cost that may capture, for instance, the cost of optic fiber to connect the sensors to the actuators, we define the cost matrix $W^K \in \mathbb{R}_+^{p \times m}$, where $W^K_{ij}$ represents the cost of establishing a communication channel from the $j$-th sensor to the $i$-th actuator. Subsequently, the last problem addressed in this paper is the minimum cost resilient co-design of \mbox{actuator-sensor-communication} described as follows.

\noindent $\mathcal P_3$  Given selective structural patterns $\bar E, \bar A \in \{0, \times,$ $ \otimes  \}^{n \times n}$,  a maximum of $k$ failures in each of actuators, sensors and communication, and  {actuation-sensing-communication} cost structure $W^B\in\mathbb{R}_+^{n \times(k+1)n}$,  $W^C\in\mathbb{R}_+^{(k+1)n \times n}$ and $W^K\in\mathbb{R}_+^{(k+1)n \times (k+1)n}$ respectively, determine $\left(\bar B^*, \bar C^*, \bar K^* \right)$ that solve  
\begin{align}
\mathop {\min }\limits_{{\begin{smallmatrix} {\bar B \in {{\left\{ {0, \times } \right\}}^{n \times \left(k+1\right)n}}}\\ {\bar C \in {{\left\{ {0, \times } \right\}}^{\left(k+1\right)n \times n}}}\\ {\bar K \in {{\left\{ {0, \times } \right\}}^{\left(k+1\right)n \times \left(k+1\right)n}}} \end{smallmatrix}}} & \qquad\left\| {\bar B} \right\|_{{W^B}}+\left\| {\bar C} \right\|_{{W^C}}+\left\| {\bar K} \right\|_{{W^K}}\\ \text{s.t.}\qquad & \hspace{-0.7cm} \left( {\bar E, \bar A,\bar B\left( {{\mathcal{I}_B}} \right), \bar C\left( {{\mathcal{I}_C}} \right), \bar K\left( {{\mathcal{K}}} \right)} \right) \text{ have no SSSFM, for all } \mathcal{I}_B, \ \mathcal{I}_C, \ \mathcal{I}_B', \ \mathcal{I}_C' \subset \mathcal{N}, \notag \\ 
&\hspace{-0.7cm}\text{ and }\mathcal{K}, \ \mathcal{K}' \subset \mathcal{N} \times \mathcal{N}, \text{with } \mathcal{I}_B = \mathcal{N} \backslash  \mathcal{I}_B',  \mathcal{I}_C = \mathcal{N} \backslash  \mathcal{I}_C', \ \  \mathcal{K} = {\left(\mathcal{N} \times \mathcal{N}\right)} \backslash  \mathcal{K}', \notag \\ 
&\hspace{-0.7cm} \text{ and } {\left| {\mathcal{I}_B'} \right| + \left| {\mathcal{I}_C'} \right| + \left| {\mathcal{K}'} \right| \le k},\notag
\end{align}
where the structural matrix $\bar K \left( {{\mathcal{K}}} \right)$ is a matrix with indices of nonzero entries contained in the set $\mathcal{K}$, and the set $\mathcal{K}'$ contains the indices of the entries representing malfunctioned communication link. 
\hfill $\diamond$

\vspace{-0.3cm}

%
%

\section{Selective Strong Structural Minimum Cost Resilient Co-Design Framework}\label{main_result}
To characterize the solutions to the problems $\mathcal{P}_1, \ \mathcal{P}_2,$  and $\mathcal{P}_3$, we proceed as follows. First, we introduce some core definitions and structures used to explicitly provide the solutions to the proposed problems. Secondly,  we provide the solutions to $\mathcal{P}_1, \ \mathcal{P}_2, \text{ and } \mathcal{P}_3$, under the assumptions that no resiliency is required, and the \mbox{actuation-sensing-communication} cost is homogenous, which we denote by $\mathcal{P}_1^0, \ \mathcal{P}_2^0, \text{ and } \mathcal{P}_3^0$ respectively, and which solutions are characterized in Theorem~\ref{lemma_sparsest}, Theorem~\ref{SSSOthm} and Theorem~\ref{theorem_SSSFM}, respectively.  Next,  we consider the problems obtained by taking into account the resiliency requirements in  $\mathcal{P}_1, \ \mathcal{P}_2, \text{ and } \mathcal{P}_3$ (under homogenous cost assumptions), which we denote by $\mathcal{P}_1^r, \ \mathcal{P}_2^r, \text{ and } \mathcal{P}_3^r$ respectively, which solutions are provided in Theorem~\ref{robustSSSC}, Corollary~\ref{robustSSSO} and Theorem~\ref{theorem_SSSFMrobust}, respectively. Lastly,  the homogenous cost assumption in $\mathcal{P}_1^r, \ \mathcal{P}_2^r, \text{ and } \mathcal{P}_3^r$ is waived, and the characterization of the general solutions to $\mathcal{P}_1, \ \mathcal{P}_2, \text{ and } \mathcal{P}_3$ is provided in Theorem~\ref{robustSSSCcost}, Corollary~\ref{robustSSSOcost} and Theorem~\ref{theorem_SSSFMcost}, respectively.

{
First, we make use of the following characterization of step differences.


\begin{definition}[Pivot and Normal Form] Given a stair matrix, a \emph{pivot} is a nonzero entry in the left-top most entry of a step difference. A step difference $\Delta_{i+1}^{i}$ of a stair matrix $M$ is \emph{normalizable} if there exist two permutation matrices $P_r^{\Delta}$ and $P_c^{\Delta}$ such that $P_r^{\Delta}\Delta_{i+1}^{i}P_c^{\Delta}$ has a pivot. Furthermore, we say that a step difference is in its \emph{normal form} if it has a pivot.
 \hfill $\diamond$
\label{pivotdef}
\end{definition}

Moreover, we can order (and label) the pivots by the induced order of the steps. Specifically, we say that two pivots $k_1$ and $k_2$ are consecutive, if there exists no other pivot $k^{\prime}$ such that $k_1<k'<k_2$.
 Furthermore, the notion of pivot will be crucial to characterize the different solutions to our problems.

Also, we require the notion of a \emph{ramp matrix}, that will enable us to characterize the feasibility space of our problems.


\begin{definition}[Ramp Structure]
A \emph{ramp structure} is a stair matrix $M\in  \{0,\times,\otimes \}^{m_1\times m_2}$ that contains a selective structural matrix $S\in  \{0,\times,\otimes \}^{n\times n}$, where $n=\min\{m_1,m_2\}$, and with $n$ step differences with pivots. \hfill $\diamond$
\label{rampmatrix}
\end{definition}

\begin{remark}\label{remark1}
From Definition~\ref{rampmatrix}, it follows that the ramp structure $M$ contains an $n \times n$ dimensional lower-triangular sub-matrix with nonzero entries in its diagonal, formed by the partially ordered columns (or rows) of $M$.
\hfill $\diamond$
\end{remark}

}

First, we provide a feasibility characterization  of $\mathcal P_1^0$, i.e., a necessary and sufficient condition to ensure SSSC.

\begin{theorem}\textbf{\emph{(Selective Strong Structural Controllability)}}
Consider the selective structural matrices $\bar E, \bar A,\bar B \in {\left\{ {0, \times, \otimes} \right\}^{n \times n}}$. The tuple $(\bar E, \bar A,\bar B)$ is SSSC  if there exist permutation matrices $P_r$ and $P_c$ such that $P_r \left [\bar A^{\lambda}\ \bar B\right] P_c$ is a ramp structure.
\hfill $\diamond$
\label{SSSCthm}
\end{theorem}

Let $\mathbb{I}_n \in \left\{ { 0, \times} \right\}^{n \times n}$ be a structural representation of an $n\times n$ identity matrix, i.e., all diagonal entries are nonzero  and the off-diagonal entries are zero. Then, every column of the structural matrix $\mathbb{I}_n$ is associated with a dedicated actuator, and  ${\mathbb{I}_n^c\left( {\mathcal I} \right) \in {\left\{ {0, \times } \right\}^{n \times \left| {\mathcal I} \right|}}}$ denotes a subset of columns in $\mathbb{I}_n$ with indices in the set $\mathcal I$. In fact, we allow the set $\mathcal I$ to be multi-index in the sense that it might contain more than once a given index. Subsequently, an input structural matrix constructed using dedicated actuators, with labels in the set $\mathcal I \subseteq \left\{ {1,2,3,\ldots, n} \right\}$, is represented by $\bar B = \left[ {\begin{array}{*{20}{c}} {\mathbb{I}_n^c\left( {\mathcal I} \right)}&{{{\bf{0}}_{n \times \left( {\left(k+1\right)n - \left| {\mathcal I} \right|} \right)}}} \end{array}} \right]$ (up to permutation), and it is referred to as \emph{dedicated solution}. Therefore, if an index in $\mathcal I$ is repeated, then it corresponds to different dedicated actuators controlling the same state variable. 

Now, consider the design objective in $\mathcal P_1^0$, under the additional restriction of dedicated actuators, i.e., $\bar B$ can have at most one non-zero entry in each column. Then, the problem reduces to that of determing the smallest set of labels $\mathcal I \subseteq \left\{ {1,2,3,\ldots, n} \right\}$ to ensure that the concatenated matrix $\left[ {\begin{array}{*{20}{c}} {{{\bar A}^\lambda }}&{\bar B} \end{array}} \right]$ can be permuted to contain a ramp structure; hence, yielding SSSC by invoking Theorem~\ref{SSSCthm}. Specifically, we obtain the following characterization of dedicated solutions to $\mathcal P_1^0$.


\begin{theorem}
Let $\bar E,\bar A\in\{0,\times,\otimes\}^{n\times n}$ be selective structural matrices, and $\bar M=P_r\bar A^{\lambda}P_c$ be a stair matrix with normalizable step differences in the normal form with $k'$ pivots, where $\bar A^{\lambda}=\bar A-\lambda \bar E$ and $\lambda \in\mathbb{C}$,  $P_r$ and $P_c$ permutation matrices with appropriate dimensions, and $\mathcal I=\{1,\ldots, n\}$ the index set. 
Then, $\bar B^*_{\Pi^r(\bar M)}=[P_r^{-1}\mathbb{I}^c_n(\mathcal I\setminus\Pi^r(\bar M)) \ \mathbf{0}_{n\times (\left(k+1\right)n-|\mathcal I\setminus\Pi^r(\bar M)|)}]$ is a dedicated solution to~$\mathcal P_1^0$, where $\Pi^r(\bar M)=\bigcup_{i=1}^{k'}p_r^i$ and $p_r^i$ denote the index of the row in $\bar M$ associated with the $i$-th pivot.
\hfill $\diamond$
\label{dedicatedSol}
\end{theorem}

Intuitively, Theorem~\ref{dedicatedSol} states that one needs to add canonical columns (to be associated with dedicated actuators) such that $\bar M$ has the columns without the step pivots are `replaced' and, subsequently, $[\bar A_S^\lambda  \quad P_r^{-1}\bar B^*_{\Pi(\bar M)}]$ can be permuted to a ramp matrix. Furthermore, multiple solutions are possible depending on the set of pivots ${\Pi(\bar M)}$ induced by $\bar M$, as we emphasize in the following remark.

\begin{remark}\label{multiplicityOfDedSol}
There are as many possible solutions as the possible combinations of pivots across different step differences. In particular, notice that a reduced number of steps differences with pivots increase the number of possible pivots for the corresponding step difference. On the other hand, if there are as many step differences as the number of rows of $\bar A^{\lambda}$, then the solution is unique and has as many dedicated actuators as the steps that do not have pivots. Lastly, observe that whereas under the homogenous cost restriction, any solution incurs in the same cost, the same is not true when such cost assumption is waived. In fact, when characterizing the solution to $\mathcal P_1,\mathcal P_2$ and $\mathcal P_3$, we will leverage this insight to consider a specific collection of pivots. \hfill $\diamond$
\end{remark}

In fact, in the next solution we show that all solutions to $\mathcal{P}_1^0$ need to be dedicated solutions.


\begin{theorem}
Given $\bar E, \bar A \in \left\{ {0, \times, \otimes} \right\}^{n \times n}$,  $\bar B^*$ is a solution to $\mathcal{P}_1^0$ if and only if $\bar B^*$ is a sparsest dedicated solution attaining SSSC.
\label{lemma_sparsest}
\hfill $\diamond$ 
\end{theorem}

In what follows, and similar to the duality between controllability and observability in linear \mbox{time-invariant} descriptor systems, one can obtain the following result.


\begin{lemma}
The tuple $\left ( \bar E, \bar A, \bar C\right)$ is SSSO if and only if $\left ( \bar E^\intercal, \bar A^\intercal, \bar C^\intercal \right)$ is SSSC.
\hfill $\diamond$
\label{duality}
\end{lemma}


By invoking Lemma~\ref{duality}, the solution to $\mathcal{P}_2^0$ can be characterized as follows.

\begin{theorem}
The structural matrix $\bar C^*$ is a solution to $\mathcal{P}_2^0$ with $\bar E, \bar A \in \left\{ {0,\times,\otimes} \right\}^{n \times n}$,  and with homogeneous sensing cost $W^C$,  if and only if $\left(\bar C^*\right)^\intercal $ is a solution to $\mathcal{P}_1^0$ with $\bar E^\intercal$ and $\bar A^\intercal$, and with homogeneous actuation cost $W^B=(W^C)^\intercal$.
\hfill $\diamond$ 
\label{SSSOthm}
\end{theorem}

Simply speaking, from Theorem~\ref{SSSOthm} and Theorem~\ref{lemma_sparsest}, it follows that the solution to $\mathcal{P}_2^0$, with homogeneous sensing cost and no sensing failures, consists of dedicated sensors, i.e., sensors that measure exactly one state variable. In fact, the dedicated solution to $\mathcal P_2^0$ can be described as follows.

\begin{corollary}
Let $\bar E,\bar A\in\{0,\times,\otimes\}^{n\times n}$ be selective structural matrices, and $\bar M=P_r\bar A^{\lambda}P_c$ be a stair matrix with normalizable step differences in the normal form with~$k'$ pivots, where $\bar A^{\lambda}=\bar A-\lambda \bar E$ and $\lambda \in\mathbb{C}$,  $P_r$ and $P_c$ permutation matrices with appropriate dimensions, and $\mathcal I=\{1,\ldots, n\}$ the index set. 
Then, $\bar B^*_{\Pi^c(\bar M)}=[\mathbb{I}^c_n(\mathcal I\setminus\Pi^c(\bar M))P_c^{-1} \ \mathbf{0}_{n\times (\left(k+1\right)n-|\mathcal I\setminus\Pi^c(\bar M)|)}]$ is a dedicated solution to  $\mathcal P_2^0$, where $\Pi^c(\bar M)=\bigcup_{i=1}^{k'}p_c^i$ and $p_c^i$ denote the index of the column in $\bar M$ associated with the $i$-th pivot.
\hfill $\diamond$
\label{dedicatedSolOutput}
\end{corollary}

Therefore, from Theorem~\ref{dedicatedSol} and Corollary~\ref{dedicatedSolOutput}, the number of dedicated actuators and sensors to ensure SSSC and SSSO, respectively, has to be the same, as formally described in the next result.


\begin{lemma}
Given solutions $\bar B^* = \left[ {\begin{array}{*{20}{c}} {\mathbb{I}_n^c\left( {\mathcal I} \right)}&{{{\bf{0}}_{n \times \left( {\left(k+1\right)n - \left| {\mathcal I} \right|} \right)}}} \end{array}} \right]$ and $\bar C^* = \left[ {\begin{array}{*{20}{c}} \left({{{\mathbb{I}_n^r} {\left(\mathcal J\right)}}}\right)^\intercal& {{{\bf{0}}_{ \left( {\left(k+1\right)n - \left| {\mathcal J} \right|} \right) \times n}^\intercal}} \end{array}} \right]^\intercal$ (both up to permutations)  to  $\mathcal{P}_1^0$ and ${\mathcal P}_2^0$, respectively, then $\left| {\mathcal I} \right| = \left| {\mathcal J} \right|$.
\hfill $\diamond$
\label{lemma_cardinality}
\end{lemma}

To obtain solution to $\mathcal{P}_3^0$, notice that the design of the information pattern will be influenced by the selection of the actuators and the sensors that are the solution to $\mathcal{P}_1^0$ and $\mathcal{P}_2^0$.  Specifically, this relationship is captured by the  following result.


\begin{proposition} Let $\mathcal{P}=\left\{ {1,2,3,\ldots,p} \right\}$ denote the labels of the actuators, $\mathcal M=\left\{ {1,2,3,\ldots,m} \right\}$ denote the labels of the sensors, and ${{\mathcal J_i}\left( {\bar K} \right)}$ contain labels of all possible sensors fed to the actuator $i$ with respect to the information pattern $\bar K$. The tuple $\left( {\bar E, \bar A, \bar B, \bar C; \bar K} \right)$ has a SSSFM (with respect to the information pattern $\bar K$) if and only if there exists {$ (E,A)\in \left(\left[\bar E\right],\left[\bar A\right]\right)^{\star}$}, $B \in \left [ \bar B \right]$, $C \in \left [ \bar C \right]$ and $K \in \left [ \bar K \right]$ and a subset of actuators $\mathcal I  \subseteq  \mathcal P$ and subset of sensors $\mathcal J \left(\mathcal I, \bar K \right)\subset \mathcal M$, described by \vspace{-0.2cm}
\begin{equation}
{\mathcal J \left(\mathcal I, \bar K \right)} = \bigcup \limits_{i \in \mathcal P \backslash \mathcal I} {{\mathcal J_i}\left( {\bar K} \right)},
\label{sen_unavail}
\end{equation}
i.e., it contains the labels of those sensors that are not fed to the actuators with labels in the set $\mathcal I$,  such that following condition holds:
\begin{equation}
rank\left[ {\begin{array}{*{20}{c}}
{A - \lambda E}&{{B \left({\mathcal I}\right)}}\\
{{C \left({{\mathcal J}\left( {{\mathcal I},\bar K} \right)}\right)}}&{\bf{0}}
\end{array}} \right] < n,  \ \lambda \in \mathbb C,
\label{strong_rank_con}
\end{equation}
where $\bf{0}$ is a zero matrix with appropriate dimensions.  \hfill $\diamond$
\label{proposition_SSSFM}
\end{proposition}

Furthermore, Proposition~\ref{proposition_SSSFM} only holds if both SSSC and SSSO are verified. In other words, SSSC and SSSO are required to attain feasibility of $\mathcal P_3^0$, as formalized in the following result.


\begin{lemma}
A system $\left( {\bar E, \bar A, \bar B, \bar C; \bar K} \right)$ has no selective strong structural fixed modes only if $\left( {\bar E, \bar A, \bar B} \right)$ is SSSC and $\left( {\bar E, \bar A, \bar C} \right)$ is SSSO.
\hfill $\diamond$
\label{lemma_SSSC_SSSO}
\end{lemma}

Hitherto, we characterized the (dedicated) solutions to $\mathcal{P}_1^0$ and $\mathcal{P}_2^0$, and also proved that SSSC and SSSO are necessary conditions for $\left( {\bar E, \bar A, \bar B, \bar C; \bar K} \right)$ to not have SSSFM. Furthermore, from Lemma~\ref{lemma_cardinality}, the number of dedicated actuators and sensors, i.e., the cardinalities of the sets $\mathcal I$ and $\mathcal J$ are the same. Besides, these are associated with a particular collection of pivots (see Theorem~\ref{dedicatedSol} and Corollary~\ref{dedicatedSolOutput}). Next, we leverage these insights to provide the pairing between sensors and actuators that ensure that $\left( {\bar E, \bar A, \bar B, \bar C; \bar K} \right)$ has no SSSFM.


{

Now, we need to introduce the notion of \emph{index-mates} that plays a key role in describing the solutions to $\mathcal P_3^0$ (as well as $\mathcal P_3^r$ and $\mathcal P_3$).

\begin{definition}[Index-mates]
Consider a stair matrix $\bar M\in\{0,\times,\otimes\}^{n\times n}$, where all the step differences are in its normal form. Also, let $\mathbb I^r_{n}(\mathcal I)$ and $\mathbb I^c_n(\mathcal J)$ be such that $\bar M_1=[\bar M \ \mathbb I^c_{n}(\mathcal I)]$ and $M_2=[\bar M^\intercal \ (\mathbb I_{n}^r(\mathcal J))^\intercal]^\intercal$ contain ramp matrices. Furthermore, let these ramp matrices be denoted by $\bar M_1(\mathcal I')$ and $\bar M_2(\mathcal J')$, where $\mathcal I'\subset \mathcal I$ and  $\mathcal J'\subset \mathcal J$, containing  the columns and rows containing the pivots  of $\bar M$, respectively, as well as the column $\mathbb{I}_n^c(\mathcal I')$ and rows $\mathbb{I}_r^c(\mathcal J')$, respectively. Then the diagonal entries corresponding to pivots are the same for $\bar M_1(\mathcal I')$ and $\bar M_2(\mathcal J')$, whereas for the remaining diagonal entries there exist an index $i\in\mathcal I'$ and $j\in\mathcal J'$ associated with the same diagonal entry, and we refer to $(i,j)$ as index-mates. \hfill $\diamond$
\label{indexMatesDef}
\end{definition}

Subsequently, let $\mathcal I =\left\{ {{i_1},{i_2},{i_3},\ldots,{i_p}} \right\}$ and $\mathcal J =\left\{ {{j_1},{j_2},{j_3},\ldots,{j_p}} \right\}$ define the indices of the effective actuators and effective sensors associated with $\bar B^*$ and $\bar C^*$ that are the solutions to $\mathcal P_1^0$ and $\mathcal P_2^0$, respectively.  Also, consider $\bar A^{\lambda}=\bar A-\lambda \bar E$ be a stair matrix where all the step differences are in its normal form. If $\left (i_{\alpha}, j_{\beta} \right)$  are index-mates when there exist  permutations $P_c$ and $P_r$ of appropriate dimensions, such that
\begin{equation}
{\left[ {{P_r}\left[ {\begin{array}{*{20}{c}}
{\bar A^{\lambda}}&{\bar B\left( {\mathcal I} \right)}\\
{\bar C\left( {\mathcal J} \right)}&{\bf{0}}
\end{array}} \right]} \right]_{ij}} = {\left[ {\left[ {\begin{array}{*{20}{c}}
{\bar A^{\lambda}}&{\bar B\left( {\mathcal I} \right)}\\
{\bar C\left( {\mathcal J} \right)}&{\bf{0}}
\end{array}} \right]{P_c}} \right]_{ij}},
\end{equation} 
where the pair of indices $\left(i,j\right)$ corresponds to the location of the diagonal entry in a step difference of the stair matrix $\bar A^{\lambda}$, with the exception of the indices of the pivots. Intuitively, we are using a dedicated actuator indexed by~$i_{\alpha}$ and dedicated sensor indexed by $j_{\beta}$ to form the ramp matrix required to ensure both SSSC and SSSO. 
}


As mentioned, the index-mates enable us to characterize the solution to $\mathcal{P}_3^0$ as follows. 

\begin{theorem}
The tuple $(\bar B^*,\bar C^*,\bar K^*)$ is a solution to $\mathcal P_3^0$ if and only if
\begin{enumerate}
\item[(\emph{i})] $\bar B^* = \left[ {\begin{array}{*{20}{c}} {\mathbb{I}_n^c\left( {\mathcal I} \right)}&{{{\bf{0}}_{n \times \left( {\left(k+1\right)n - \left| {\mathcal I} \right|} \right)}}} \end{array}} \right]$ is a solution to $\mathcal{P}_1^0$ and $\bar C^* = \left[ {\begin{array}{*{20}{c}} {{\left(\mathbb{I}_n^r{\left(\mathcal J\right)}\right)}^\intercal}& {{{\bf{0}}_{ \left( {\left(k+1\right)n - \left| {\mathcal J} \right|} \right) \times n}^\intercal}} \end{array}} \right]^\intercal$ is a solution to ${\mathcal P}_2^0$ (both up to permutations);
\vspace{0.2cm}
\item[(\emph{ii})] $\bar K^*_{i_{\alpha}, j_{\beta}}=\times$, if the pair $\left (i_{\alpha}, j_{\beta} \right)\in\mathcal I'\times \mathcal J'$ are index-mates, and zero otherwise, where $\mathcal I' =\left\{ {{i_1},{i_2},{i_3},\ldots,{i_p}} \right\}$ and  $\mathcal J' =\left\{ {{j_1},{j_2},{j_3},\ldots,{j_p}} \right\}$ are the indices of the effective actuators associated with $\bar B^*$ and effective sensors associated with $\bar C^*$, respectively. \hfill $\diamond$
\end{enumerate}
 \label{theorem_SSSFM}
\end{theorem}

%

To achieve robustness to $k$ actuator-sensor-communicati- on failures, a possible approach to achieve feasibility for $\mathcal{P}_1^r, \ \mathcal{P}_2^r, \text{ and } \mathcal{P}_3^r$ is that of considering $k+1$ sparsest solutions to $\mathcal{P}_1^0, \ \mathcal{P}_2^0, \text{ and } \mathcal{P}_3^0$. Notably, such strategy results in optimal solutions in strong structural theory, which contrast to the solutions to similar problems under the requirement of (not strong) structural systems theory conditions (see, for instance,~\cite{Liu,2015arXiv150707205L}). 


\begin{theorem}
Let $\{\bar B_i^*=[\mathbb{I}_n^c(\mathcal I_i) \ {\mathbf{0}}_{n \times \left(\left(k+1\right)n-\left| \mathcal I_i \right| \right)}]\}_{i=1}^{k+1}$ be a collection of solutions to $\mathcal P_1^0$. The solution to~$\mathcal P_1^r$ is given by $\bar B^*=[\mathbb{I}_n^c(\bigcup_{i=1}^{k+1}\mathcal I_i) \ \mathbf{0}_{n \times \left(\left(k+1\right)n-\# \right)}]$ (up to permutation of the columns), where $\#={\sum\limits_{i = 1}^{k + 1}{\left| \mathcal I_i \right|}}$.~\hfill~$\diamond$
\label{robustSSSC}
\end{theorem}

Once again, invoking duality (see Theorem~\ref{SSSOthm}), we obtain the following result.

\begin{corollary}
Let $\{\bar C_i^*=[{\left(\mathbb{I}_n^r(\mathcal J_i)\right)}^\intercal \ \mathbf{0}_{\left(\left(k+1\right)n-\left| \mathcal J_i \right| \right)\times n}^\intercal]^\intercal\}_{i=1}^{k+1}$ be a collection of solutions to $\mathcal P_2^0$. The solution to $\mathcal P_2^r$ is given by $\bar C^*=[{\left(\mathbb{I}_n^r(\bigcup_{i=1}^{k+1}\mathcal J_i)\right)}^\intercal  \ \mathbf{0}_{\left(\left(k+1\right)n-\# \right)\times n}^\intercal]^\intercal$ (up to permutation of the rows), where $\#={\sum\limits_{i = 1}^{k + 1}{\left| \mathcal J_i \right|}}$.~\hfill~$\diamond$
\label{robustSSSO}
\end{corollary}

Finally, the solution to $\mathcal P_3^r$ can be described as follows.


\begin{theorem}
The tuple $(\bar B^*,\bar C^*,\bar K^*)$ is a solution to $\mathcal P_3^r$  if and only if 
\begin{enumerate}
\item[(\emph{i})] $\bar B^*$ and $\bar C^*$ are solutions to $\mathcal{P}_1^r$ and ${\mathcal P}_2^r$, respectively; and 
\item[(\emph{ii})]$\bar K^*_{i_{\alpha}, j_{\beta}}=\times$, when $\left (i_{\alpha}, j_{\beta} \right)\in\mathcal I'\times \mathcal J'$ are index-mates, where $\mathcal I' =\left\{ {{i_1},{i_2},{i_3},\ldots,{i_{p'}}} \right\}$ and  $\mathcal J' =\left\{ {{j_1},{j_2},{j_3},\ldots,{j_{p'}}} \right\}$ are the indices of the effective actuators and the effective sensors, respectively, and zero otherwise.\hfill $\diamond$
\end{enumerate}
\label{theorem_SSSFMrobust}
\end{theorem}
%

Now, recall Remark~\ref{multiplicityOfDedSol} where we emphasized that specific selection of pivots in the normalized stair matrix  leads to different dedicated solutions, and overall actuation cost. In particular, following the different constructions provided up to this point,  the dedicated actuators required complement the existence of the pivots to form a ramp matrix. Therefore, suppose that each entry $[\bar A^\lambda]_{i,j}$ is associated with cost $[ W^B]_{i,j}$, then we can say that a pivot has the cost associated with an entry in the actuation cost matrix. Hence, we can associate a collection of pivots with an overall cost. Therefore, it follows that the selected pivots should incur in the maximum actuation cost, which implies that the collection of dedicated actuators will incur in the minimum cumulative actuation cost.  We formalize these observations in the following result, where we provide the characterization of the solutions to $\mathcal P_1$.
 

\begin{theorem}
Let $\{\bar B_i^*=[\mathbb{I}_n^c(\mathcal I_i) \ {\mathbf{0}}_{n \times \left(\left(k+1\right)n-\left| \mathcal I_i \right| \right)}]\}_{i=1}^{k+1}$ be a collection of solutions to $\mathcal P_1^0$ (up to permutation of the columns), constructed as in Theorem~\ref{dedicatedSol}, where the sum of the pivots' cost described by $W^B$ is maximized. Then $\bar B^*=[\mathbb{I}_n^c(\bigcup_{i=1}^{k+1}\mathcal I_i) \ \mathbf{0}_{n \times \left(\left(k+1\right)n-\# \right)}]$ is a solution to $\mathcal P_1$ (up to permutation of the columns), where $\#={\sum\limits_{i = 1}^{k + 1}{\left| \mathcal I_i \right|}}$.
~\hfill~$\diamond$
\label{robustSSSCcost}
\end{theorem}

Similarly, invoking duality and Theorem~\ref{robustSSSCcost}, and associating the cost of the pivots with the sensing cost matrix, we obtain the following result.

\begin{corollary}
Let $\{\bar C_i^*=[{\left(\mathbb{I}_n^r(\mathcal J_i)\right)}^\intercal \ \mathbf{0}_{\left(\left(k+1\right)n-\left| \mathcal J_i \right| \right)\times n}^\intercal]^\intercal\}_{i=1}^{k+1}$ be a collection of solutions to $\mathcal P_2^0$ (up to permutation of the rows), which solutions were computed using duality and Theorem~\ref{dedicatedSol}, where the sum of the pivots' cost described by $W^C$ is maximized. Then  $\bar C^*=[{\left(\mathbb{I}_n^r(\bigcup_{i=1}^{k+1}\mathcal J_i)\right)}^\intercal  \ \mathbf{0}_{\left(\left(k+1\right)n-\# \right)\times n}^\intercal]^\intercal$ is a solution to $\mathcal P_2$ (up to permutation of the rows), where $\#={\sum\limits_{i = 1}^{k + 1}{\left| \mathcal J_i \right|}}$.~\hfill~$\diamond$
\label{robustSSSOcost}
\end{corollary}


In order to construct the solution to $\mathcal P_3$, we consider the solutions to $\mathcal P_3^0$ and determine  one that incurs in the minimum \mbox{actuation-sensing-communication} cost. Now, notice that  since the solutions to $\mathcal P_3^0$ establish a one-to-one correspondence between dedicated sensors and actuators through a single communication channel, one can associate a cost with index-mates $(i,j)$ that comprises the cost of actuating the state variable controlled by the dedicated actuator $i$, the cost of sensing the state variable measured by the dedicated sensor $j$, and the communication from dedicated sensor $j$ to dedicated actuator $j$. Therefore, the different index-mates have a cost associated with  $(W^B,W^C,W^K)$, and it follows that the solution to~$\mathcal P_3$ consists of  determining the solutions to $\mathcal P_3^0$ that incur minimum \mbox{actuation-sensing-communication} (or, equivalently, the index-mates) cost. Specifically, the minimum cost solutions to $\mathcal P_3$ are constructed as we formally state in the following theorem.

\begin{theorem}
Let $\bar B^* $ and $\bar C^*$ be the solutions to $\mathcal{P}_1^r$ and ${\mathcal P}_2^r$, respectively, obtained by permuting $\bar A^\lambda$ to a normal form with largest sum of pivots' cost  described by $W^B+W^C+W^K$. Then the tuple $(\bar B^*,\bar C^*,\bar K^*)$ is a solution to $\mathcal P_3$, where $\bar K^*_{i_{\alpha}, j_{\beta}}=\times$ for the collection of index-mates $\left (i_{\alpha}, j_{\beta} \right)\in\mathcal I'\times \mathcal J'$  that incur the minimum overall cost  associated with $(W^B,W^C,W^K)$, with $\mathcal I' =\left\{ {{i_1},{i_2},{i_3},\ldots,{i_{p'}}} \right\}$ and $\mathcal J' =\left\{ {{j_1},{j_2},{j_3},\ldots,{j_{p'}}} \right\}$ denoting the indices of the effective actuators and effective sensors respectively, and zero otherwise. 
\hfill $\diamond$
\label{theorem_SSSFMcost}
\end{theorem}

\vspace{-0.3cm}

\section{Illustrative Example}\label{illus_examp}
{In this section, we illustrate the use of the main results regarding  the minimum cost resilient actuation-sensing-communication \mbox{co-design} in the context of the electric power grid. Specifically, we consider a $5$-bus example in~\cite{camsap}, whose dynamics is approximated by a $16$-th order multi-input multi-output linear time invariant system (i.e., a particular case of descriptor systems), and which state variables description can be found in Table~I in~\cite{camsap}.

A normalized stair matrix $\bar M_{5_\text{bus}}^{\lambda}$ is given as follows 
\begin{align*}
{}    &     {\begin{smallmatrix} & & 3 & {\color{blue}2} & & 1&{12}&{6}&{\color{blue}5}&{4}&{13}&{9}&{\color{blue}8}&{7}&{10}&{14}&{15}&{11}&{16} \end{smallmatrix}}\\
{\begin{smallmatrix} {2}\\ {3}\\ {1}\\ {5}\\ {6}\\ {4}\\ {8}\\ {9}\\ {12}\\ {13}\\ {7}\\ {10}\\ {14}\\ {15}\\ {11}\\ {16} \end{smallmatrix}}     &     
{\left[ {\begin{smallmatrix}
 \times & \otimes & {}&{}&{}&{}&{}&{}&{}&{}&{}&{}&{}&{}&{}&{}\\
 \otimes &0& \times &{}&{}&{}&{}&{}&{}&{}&{}&{}&{}&{}&{}&{}\\
 \times & \times & \otimes & \times &{}&{}&{}&{}&{}&{}&{}&{}&{}&{}&{}&{}\\
0&0&0&0& \times & \otimes &{}&{}&{}&{}&{}&{}&{}&{}&{}&{}\\
0&0&0&0& \otimes &0& \times &{}&{}&{}&{}&{}&{}&{}&{}&{}\\
0&0&0&0& \times & \times & \otimes & \times &{}&{}&{}&{}&{}&{}&{}&{}\\
0&0&0&0&0&0&0&0& \times & \otimes &{}&{}&{}&{}&{}&{}\\
0&0&0&0&0&0&0&0& \otimes &0& \times &{}&{}&{}&{}&{}\\
0&0& \times & \otimes &0&0& \times &0&0&0&0& \times &{}&{}&{}&{}\\
0&0& \times &0&0&0& \times & \otimes &0&0& \times & \times &{}&{}&{}&{}\\
0&0&0&0&0&0&0&0&0&0& \otimes &0& \times &{}&{}&{}\\
0&0&0&0&0&0&0&0&0&0&0& \otimes &0& {\ \times \ } &{}&{}\\
0&0&0&0&0&0& \times &0&0&0& \times &0& \otimes &0& \times &{}\\
0&0& \times &0&0&0& \times &0&0&0&0& \times &0& \otimes & \times &{}\\
0&0&0&0&0&0&0&0&0&0&0&0&0&0& \otimes & {\ \times \ } \\
0&0&0&0&0&0&0&0&0&0& \times & \times &0&0& \times & \otimes 
\end{smallmatrix}} \right],}\vspace{-0.5cm}
\end{align*}
where the row and the column indices labelled in $\bar M_{5_\text{bus}}^{\lambda}$  correspond to the  row and the column indices in  $\bar A_{5_\text{bus}}^{\lambda}$ before permutations, respectively. Simply speaking, they correspond to the state variables indices to be considered for the actuation-sensing-communicaton design. 

Notice that $\bar M_{5_\text{bus}}^{\lambda}$ contains $13$ pivots. Therefore, as prescribed by Theorem~\ref{dedicatedSol}, it requires three dedicated actuators, and the possible collection of state variables that can be controlled by dedicated actuators are indexed by $\mathcal I_1=\{12,14,16\}$, $\mathcal I_2=\{13,14,16\}$, $\mathcal I_3=\{12,15,16\}$, or $\mathcal I_4=\{13,15,16\}$ (notice the row indices highlighted in red). As stated in Theorem~\ref{lemma_sparsest}, these dedicated solutions are also the solution to $\mathcal P_1^0$. In contrast, the solution to the $\mathcal P_2^0$ is unique, and consists in measuring by three dedicated sensors the variables indexed by $\mathcal J=\{2,5,8\}$.

Subsequently, suppose that the actuation cost associated with the possible actuation schemes in monotonic, i.e., actuating state variable $i$ is less than that of actuating state variable $j$ for $j>i$, whereas the measuring cost is finite and different for all state variables. Furthermore, assume that we want to ensure that the system is robust with respect to one failure, that is $k=1$ in $\mathcal P_1$, $\mathcal P_2$, and~$\mathcal P_3$.

Subsequently, by invoking Theorem~\ref{robustSSSC}, it follows that the largest combination of pivots leads to the set of actuated state variables described by the indices in $\mathcal I_1$, so these should be actuated twice. By invoking Corollary~\ref{robustSSSO}, due to the uniqueness of solution to $\mathcal P_2^0$, it follows that the state variables described by $\mathcal J$ need to be measured twice. Therefore, we construct the matrices $\bar B_{5_\text{bus}}=[\mathbb{I}_{16}^c(\mathcal I) \ \mathbf{0}_{16 \times 26}]$ and $\bar C_{5_\text{bus}}=[\left(\mathbb{I}_{16}^r(\mathcal J)\right)^\intercal \ \mathbf{0}_{26 \times 16}^\intercal]^\intercal$, where $\mathcal I=\{12,14,16,12,14,16\}$ and $\mathcal J=\{2,5,8,2,5,8\}$. Furthermore, the indices of effective actuators and sensors is given by $\mathcal I'=\{1,\ldots, 6\}$ and $\mathcal J'=\{1,\ldots, 6\}$, respectively.  Finally, by invoking Theorem~\ref{theorem_SSSFMcost}, we obtain  the information pattern $\bar K_{5_\text{bus}} \in \{0,\times\}^{32\times32}$ has nonzero entries corresponding to the index-mates $\{\left(1, 2\right),\left(4, 5\right),\left(2,3\right),\left(5,6\right),\left(3,1\right),\left(6,4\right)\}$ that incur minimum \mbox{actuation-sensing-communication} cost. Then, $\left(\bar B_{5_\text{bus}},\bar C_{5_\text{bus}}, \bar K_{5_\text{bus}}\right)$ is a solution to $\mathcal P_3$.

\vspace{-0.2cm}

\section{Conclusions and Further Research}\label{concl_furth}
In this paper, we introduced the selective strong structural systems framework and used it to address the problem of minimum cost resilient actuation-sensing-communication co-design for descriptor systems. We argue that such setup is ideal for coping with several scenarios where uncertainty in the system's modeling is present, and guarantees are required for any possible scenario. Furthermore, we introduced the notion of selective strong structural fixed modes as a characterization of the feasibility of decentralized control laws and provided necessary and sufficient conditions for this property to hold. Also, we showed how these conditions could be leveraged to determine the minimum cost resilient placement of actuation-sensing-communication technology such that decentralized control through static output feedback is possible, and unveiled the connection with closely related problems of minimum cost resilient actuation and sensing placement while achieving selective strong structural controllability and observability, respectively. 

Future research will address the design when diverse cost structures are considered to be associated with collection of actuators/sensors that can actuate/measure several state variables at the same time, i.e., the possible structure of the input and output matrices is restricted. Also, we aim to exploit the proposed problems' structure to derive efficient algorithms to attain the proposed designs. 

\vspace{-0.2cm}

\section*{Appendix}
%
%
\textbf{Proof of Theorem~\ref{SSSCthm}}: Assume that for the selective structural matrices $\bar E, \bar A,\bar B \in {\left\{ {0, \times, \otimes} \right\}^{n \times n}}$, there exist permutation matrices $P_r$ and $P_c$ such that $P_r \left [\bar A^{\lambda}\ \bar B\right] P_c$ is a ramp structure. Then, it contains $n \times n$ dimensional lower-triangular sub-matrix with non-zero entries in its diagonal, which implies that $
\text{rank}\left[ {\begin{array}{*{20}{c}}
{A - \lambda E}&{B}
\end{array}} \right] = n,  \ \lambda \in \mathbb C. $
Thus, by invoking the controllability criteria for regular descriptor systems (see Theorem 7 in~\cite{Yip}), it follows that $\left(E, A, B\right)$ is controllable; since, the above holds for all numerical realizations of $\left(\bar E, \bar A, \bar B\right)$, we conclude that $\left(\bar E, \bar A, \bar B\right)$, is SSSC. 
%
\hfill $\blacksquare$\\
%
%
\textbf{Proof of Theorem~\ref{dedicatedSol}}: For $\bar E,\bar A\in\{0,\times,\otimes\}^{n\times n}$, let $\bar M=P_r\bar A^{\lambda}P_c$ be a stair matrix with step differences in the normal form, where $\bar A^{\lambda}=\bar A-\lambda \bar E$, $\lambda \in\mathbb{C}$,  $P_r$ and $P_c$ are the permutation matrices with appropriate dimensions. Furthermore, let $p_r^i$ and $p^{i}_c$ denote the row and the column indices of the $i$-th pivot in $\bar M$, respectively, with $i=1,\ldots, k'$. In addition, let $\bar M_{p^{i}_c}$ represent the column in $\bar M$ containing $i$-th pivot. Consider $\bar B^*_{\Pi^r(\bar M)}=[P_r^{-1}\mathbb{I}^c_n(\mathcal I\setminus\Pi^r(\bar M)) \ \mathbf{0}_{n\times (\left(k+1\right)n-|\mathcal I\setminus\Pi^r(\bar M)|)}]$, then there exists $P'_c$, such that
\begin{align}
\bar M P'_c =& \left[\bar M_{p^{1}_c} {\mathbb I_n^c}\left({{\mathcal J}\left(p^{1}_c, p^{2}_c\right)}\right) \bar M_{p^{2}_c} \ \ \bar M_{p^{2}_c} {\mathbb I_n^c}{\left({\mathcal J}\left(p^{2}_c, p^{3}_c\right)\right)} \bar M_{p^{3}_c}  \right. \nonumber \\
& \qquad \ldots \left. \bar M_{p^{k'-1}_c} {\mathbb I_n^c}{\left({\mathcal J}\left(p^{k'-1}_c, p^{k'}_c\right)\right)} \bar M_{p^{k'}_c}  \ \ \bar M^{-}\ \ \mathbf{0}_{n\times (\left(k+1\right)n-|\mathcal I|)} \right], \nonumber
\end{align}
where $\mathcal J\left(p^{i-1}_c, p^{i}_c\right)=\left\{j \in \mathbb N \ :\ p^{i-1}_c<j<p^{i}_c\right\}$ and $\bar M^{-}$ denotes the columns of $\bar M$ without pivots. Then the matrix $\left[\bar A - \lambda \bar E\ \bar B\right]$ can be permuted to a ramp structure~$\bar M$, and by invoking Theorem~\ref{SSSCthm}, $\left(\bar E, \bar A, \bar B \right)$ is SSSC. Finally, $\bar B$ contains minimum number of nonzero columns with unique row indices in $\mathcal I\setminus\Pi^r(\bar M)$ necessary to ensure the ramp structure of $\bar M$ { as consequence of Assumption~1}, which implies that $\bar B$ is a sparsest dedicated solution to $\mathcal P_1^0$.  \hfill $\blacksquare$\\
%
%
\textbf{Proof of Theorem~\ref{lemma_sparsest}}: Suppose $\bar B$ is a solution to $\mathcal{P}_1^0$ and assume that $\bar B$ is not a sparsest dedicated solution. Then, it follows that $\bar B$ contains at least one column $\bar B^i$ with more than one non-zero entry. Since $\left(\bar E, \bar A, \bar B\right)$ is SSSC, then by Theorem~\ref{SSSCthm} the matrix $\left[ {{{{\bar A}-\lambda \bar E }}\ {\bar B}} \right]$ can be permuted to a ramp matrix $\bar M$, for every $\lambda \in \mathbb C$. The ramp structure property is preserved if each non-zero entry below the pivot in $\bar B^i$ (with its entries reordered in $\bar M$) is replaced with a zero. Therefore, the column $\bar B^i$ can be replaced with a sparsest column $\bar B^{'i}$ with exactly one non-zero entry, resulting in selective structural matrix $\bar B'$ such that the tuple $\left(\bar E, \bar A, \bar B'\right)$ is SSSC. It follows that $\bar B'$ being sparser than $\bar B$ incurs lower actuation cost (under homogenous actuation cost assumption), which contradicts the hypothesis that $\bar B$ is a solution to $\mathcal{P}_1^0$. Hence, $\bar B$ is a feasible solution to $\mathcal{P}_1^0$ if it is also the sparsest dedicated solution for the tuple $\left(\bar E, \bar A, \bar B\right)$ to be SSSC. Now assume that $\bar B$ is a sparsest dedicated solution for the tuple $\left(\bar E, \bar A, \bar B\right)$ to be SSSC. Then, by invoking Theorem~\ref{dedicatedSol}, $\bar B$ is a solution to $\mathcal P_1^0$. \hfill $\blacksquare$\\
%
%
\textbf{Proof of Lemma~\ref{duality}}: Suppose the tuple $\left ( \bar E^\intercal, \bar A^\intercal, \bar C^\intercal \right)$ is SSSC. Then, by Definition~\ref{def_SSSC} and the controllability criteria for { regular} descriptor systems (Theorem 7 in~\cite{Yip}) the ${rank} \left[ {A^\intercal - \lambda E^\intercal} \ C^\intercal \right]=n$ for all {$ (E,A)\in \left(\left[\bar E\right],\left[\bar A\right]\right)^{\star}$}, $C \in \left[\bar C \right]$ and $\lambda \in \mathbb{C}$. This is equivalent to $rank\left[ {A^\intercal - \lambda E^\intercal} \ C^\intercal \right]^\intercal = n$, so by invoking Definiton~\ref{def_SSSO}, it follows that $\left ( \bar E, \bar A, \bar C\right)$ is SSSO. \hfill $\blacksquare$\\
\textbf{Proof of Proposition~\ref{proposition_SSSFM}}: Consider a numerical realization $\left( { E,  A,  B,  C; \bar K} \right)$, where {$ (E,A)\in \left(\left[\bar E\right],\left[\bar A\right]\right)^{\star}$}, $B \in \left[\bar B\right]$ and $C \in \left[\bar C\right]$. Then, from~\cite{Anderson1981703}, $\left( { E,  A,  B,  C; \bar K} \right)$ has a fixed mode with respect to (w.r.t.) the information pattern $\bar K$, if and only if there exist ${\mathcal I} \subseteq \mathcal P$ and $\mathcal J \left(\mathcal I, \bar K \right) \subseteq \mathcal M$, such that the following holds
\vspace{-0.3cm}
\begin{equation*}
rank\left[ {\begin{array}{*{20}{c}}
{A - \lambda E}&{{B \left({\mathcal I}\right)}}\\
{{C \left({{\mathcal J}\left( {{\mathcal I},\bar K} \right)}\right)}}&{\bf{0}}
\end{array}} \right] < n. \vspace{-0.3cm}
\end{equation*}
This criterion is equivalent to that presented in Definition~\ref{struct_fixedmodes}, and it implies that $\left( {\bar E,\bar  A,\bar  B,\bar  C; \bar K} \right)$ has a SSSFM w.r.t. the information pattern $\bar K$. \hfill $\blacksquare$\\
%
%
\textbf{Proof of Lemma~\ref{lemma_SSSC_SSSO}}: Let the sets $\mathcal I$ and $\mathcal J$ contain the labels of the effective actuators in $\bar B$ and the effective sensors in $\bar C$, respectively. Now, consider a scenario where the information pattern matrix $\bar K$ is full, i.e., the sensor measurements are available to all the actuators. By invoking Proposition~\ref{proposition_SSSFM} for the case when the none of the actuators with labels in $\mathcal I$ are fed the measurements from the sensor with labels in $\mathcal J$, it follows that $\mathcal J \left(\mathcal I, \bar K \right)= \mathcal J$ by~\eqref{sen_unavail}. Therefore, the rank condition in~\eqref{strong_rank_con} is the same as the observability criterion, and by invoking Definition~\ref{def_SSSO} it follows that $\left( {\bar E, \bar A, \bar B, \bar C; \bar K} \right)$ does not have selective strong structural fixed modes only if $\left( {\bar E, \bar A, \bar C} \right)$ is SSSO. Similarly, for a scenario where the information pattern matrix $\bar K$ is full, consider a case when measurements from all the sensors with labels in $\mathcal J$ are fed to the actuators with labels in $\mathcal I$, implying that $\mathcal J \left(\mathcal I, \bar K \right)=\emptyset$ {by}~\eqref{sen_unavail}. Then \eqref{strong_rank_con} results in the controllability criteria, and by invoking Definition~\ref{def_SSSC},  it follows that $\left( {\bar E, \bar A, \bar B, \bar C; \bar K} \right)$ does not have SSSFM only if $\left( {\bar E, \bar A, \bar B} \right)$ is SSSC.\hfill $\blacksquare$\\
%
%
\textbf{Proof of Theorem~\ref{theorem_SSSFM}}: Suppose the system $\left( {\bar B, \bar C, \bar K} \right)$ is a solution to $\mathcal P_3^0$, implying that it incurs minimal \mbox{actuation-sensing-communication} cost (under homogeneous cost assumption) for the system $\left( {\bar E, \bar A, \bar B, \bar C; \bar K} \right)$ to not have SSSFM, where $\bar E, \bar A \in \left\{ {0, \times, \otimes} \right\}^{n \times n}$. Then, by Lemma~\ref{lemma_SSSC_SSSO}, this is possible only if the tuples $\left( {\bar E, \bar A, \bar B} \right)$ and $\left( {\bar E, \bar A, \bar C} \right)$ are SSSC and SSSO, respectively. Furthermore, let the stair matrix $\bar M=P_r\left[\bar A-\lambda \bar E\right]P_c$ contain step differences in normal form with $k'$ pivots, where $P_r$, $P_c$ are the permutation matrices with appropriate dimensions. By invoking Theorem~\ref{dedicatedSol}, $\bar B$ must contain at least $p=n-k'$ nonzero columns (with exactly one nonzero entry in each) for $\left( {\bar E, \bar A, \bar B} \right)$ to be SSSC, i.e., $\bar B= \left[ {\begin{array}{*{20}{c}} {\mathbb{I}_n^c\left( {\mathcal I} \right)}&{{{\bf{0}}_{n \times \left( { \left(k+1\right)n - \left| {\mathcal I} \right|} \right)}}} \end{array}} \right]$ where $\left|\mathcal I\right|=p$. Hence, $\bar B$ is also a solution to $\mathcal P_1^0$ (by Theorem~\ref{lemma_sparsest}). Similarly by invoking Lemma~\ref{lemma_cardinality} and Theorem~\ref{SSSOthm}, $\bar C = \left[ {\begin{array}{*{20}{c}} {\left({\mathbb{I}_n^r}{\left(\mathcal J\right)}\right)^\intercal}& {{{\bf{0}}_{ \left( { \left(k+1\right)n - \left| {\mathcal J} \right|} \right) \times n}^\intercal}} \end{array}} \right]^\intercal$, is a solution to ${\mathcal P}_2^0$ (up to permutation of rows), where $\left|\mathcal J\right|=p$. 

Next, we show that $\|\bar K\|_0=p$. Assume without loss of generality ${\left\| {\bar K} \right\|_0}=p-1$, implying that there are $p-1$ communication channels between the sensors and the actuators. By invoking Proposition~\ref{proposition_SSSFM} w.r.t. the information pattern $\bar K$, consider the case when none of the actuators with labels in $\mathcal I$ are fed the measurements from the sensors with labels in $\mathcal J$. Then, by (\ref{sen_unavail}), $\mathcal J \left(\mathcal I, \bar K \right)=\mathcal J\setminus\{j_{\beta}\}$, and $j_{\beta} \in \mathcal J$ denotes the label of the sensor that remains non utilized, that always exist since only $p-1$ communication channels exist. Thus, there exist {$ (E,A)\in \left(\left[\bar E\right],\left[\bar A\right]\right)^{\star}$}, $C \in \left [ \bar C\left(\mathcal J \left(\mathcal I, \bar K \right)\right) \right]$, such that the condition in~\eqref{strong_rank_con} holds, i.e., $rank\left[ {A^{\intercal} - \lambda E}^{\intercal} \ {C^{\intercal}}  \right]^{\intercal} < n$ where $\lambda \in \mathbb C$. This implies that the tuple $\left( {\bar E, \bar A, \bar B, \bar C; \bar K} \right)$ has SSSFM and $\left( {\bar B, \bar C, \bar K} \right)$ is not a solution to $\mathcal P_3^0$, which is a contradiction. Hence, the condition ${\left\| {\bar K} \right\|_0} = p$ is necessary for feasibility, and it is also optimal, since we are using minimum number of dedicated actuators and dedicated sensors. This implies that  to feed the measurements from $p$ dedicated sensors to $p$ dedicated actuators, where $\bar K_{i_{\alpha}, j_{\beta}}=\times$ for the index-mates $\left (i_{\alpha}, j_{\beta} \right)$, where $\left (i_{\alpha}, j_{\beta} \right)\in\mathcal I\times \mathcal J$. Now, suppose there exist a sensor with index $j_{\beta} \in \mathcal J$ that does not have an index mate in the set $\mathcal I$. Then, by Proposition~\ref{proposition_SSSFM}, consider the scenario when $\mathcal I=\emptyset$, which implies $\mathcal J \left(\mathcal I, \bar K \right)=\mathcal J\setminus\{j_{\beta}\}$ by~\eqref{sen_unavail}, and it follows that the condition in~\eqref{strong_rank_con} holds. This implies that $\left(\bar E, \bar A, \bar B, \bar C; \bar K \right)$ has a SSSFM, which contradicts the hypothesis. Hence, corresponding to a sensor with index~$j_{\beta}$, there always exists an actuator with index $i_{\alpha} \in \mathcal I$, such that the pair $\left (i_{\alpha}, j_{\beta} \right)$ are index-mates.

Finally, we notice that the reverse implication immediately holds by reusing the same arguments presented about regarding the feasibility and optimality.
\hfill $\blacksquare$\\
%
%
\textbf{Proof of Theorem~\ref{robustSSSC}}: Assume that $\bar B=[\mathbb{I}_n^c(\bigcup_{i=1}^{k+1}\mathcal I_i) $ $ \mathbf{0}_{n \times \left(\left(k+1\right)n- \# \right)}]$ is a solution to $\mathcal P_1^r$, where $\{\bar B_i=[\mathbb{I}_n^c(\mathcal I_i) \ {\mathbf{0}}_{n \times \left(\left(k+1\right)n-\left| \mathcal I_i \right| \right)}]\}_{i=1}^{k+1}$ is a collection of solutions to $\mathcal P_1^0$ and $\#={\sum\limits_{i = 1}^{k + 1}{\left| \mathcal I_i \right|}}$. By Theorem~\ref{lemma_sparsest}, it follows that $\bar B_i$ are sparsest dedicated solution for the tuple $\left(\bar E, \bar A, \bar B_i\right)$ to be SSSC, and therefore, $\left|\mathcal I_1\right|=\left|\mathcal I_2\right| \ldots \left|\mathcal I_{k+1}\right|=p$ (by Theorem~\ref{dedicatedSol}). The nonzero columns in $\bar B_i$ (with exactly one nonzero entry each) provide $p$ pivots so that the matrix $\left[\bar A-\lambda \bar E \ \bar B_i\right]$ can be permuted to a ramp structure $\bar M_i$, where $\lambda \in\mathbb{C}$. Let the dedicated actuator be represented by the nonzero column $\bar B_i^j$, associated with the $j$-th pivot, where $j\in\left\{1,\ldots,p\right\}$. In order to permute $\left[\bar A-\lambda \bar E \ \bar B\right]$ to a ramp structure $\bar M$, $\bar B$ contains a set of $k+1$ dedicated actuators $\{\bar B_i^{j}\}_{i=1}^{k+1}$. In other words, the dedicated actuators $\{\bar B_i^{j}\}_{i=1}^{k+1}$ provide $k+1$ pivots (with the same row indices) in $\bar M$. Without loss of generality, consider a case of $k$ actuator failures represented by the set of columns $\{\bar B_i^{j}\}_{i=1}^{k}$. Then the functioning actuator $\bar B_{k+1}^{j}$ in $\bar B$, representing the column associated with the $j$-th pivot, preserves the ramp structure of $\bar M$, thereby ensuring the SSSC property of the tuple $\left(\bar E, \bar A, \bar B \right)$ after $k$ actuator failures (by Theorem~\ref{SSSCthm}). Hence, $\bar B=[\mathbb{I}_n^c(\bigcup_{i=1}^{k+1}\mathcal I_i) \ \mathbf{0}_{n \times \left(\left(k+1\right)n-\# \right)}]$ is necessary to ensure feasibility. In fact, we need a minimum of $k+1$ nonzero columns for each of the column in $\bar M$ without a pivot, to ensure robustness w.r.t. $k$ actuator failures. Hence, the solution $\bar B$ is also optimal. \hfill $\blacksquare$\\
%
%
\textbf{Proof of Theorem~\ref{theorem_SSSFMrobust}}: Assume the tuple $(\bar B,\bar C,\bar K)$ is a solution to $\mathcal P_3^r$, implying that it incurs minimal \mbox{actuation-sensing-communication} cost (under homogeneous cost assumption) for the system $\left( {\bar E, \bar A, \bar B, \bar C; \bar K} \right)$ to not have SSSFM, where $\bar E, \bar A \in \left\{ {0, \times, \otimes} \right\}^{n \times n}$. Furthermore, this property is robust w.r.t. a total of $k$ failed actuators, sensors and communication links. Consider a case of $k$ actuator failures. Then, by Lemma~\ref{lemma_SSSC_SSSO}, the SSSC property of $\left(\bar E, \bar A, \bar B \right)$ must hold under $k$ failed actuators for $\left( {\bar E, \bar A, \bar B, \bar C; \bar K} \right)$ to not have SSSFM. Hence, by similar reasoning provided in the proof of Theorem~\ref{robustSSSC}, $\bar B$ is a solution to $\mathcal P_1^r$. Similarly, by Corollary~\ref{robustSSSO}, it follows that $\bar C$ is a solution to $\mathcal P_2^r$. Furthermore, there exist $k+1$ pairs of index-mates $\left (i_{\alpha}, j_{\beta} \right)$, such that $i_{\alpha}$ and $i_{\beta}$ denote the indices of the 
effective dedicated actuator and sensor, that enable the matrices $\left[\bar A-\lambda \bar E \ \bar B \right]$ (Theorem~\ref{dedicatedSol}) and $\left[\bar A^\intercal-\lambda \bar E^\intercal \ \bar C^\intercal \right]^\intercal$ (Theorem~\ref{SSSOthm}), to be permutable to ramp structures, as described after Definition~\ref{indexMatesDef}. Now consider the information pattern $\bar K_{i_{\alpha}, j_{\beta}}=\times$, when $\left (i_{\alpha}, j_{\beta} \right)\in\mathcal I\times \mathcal J$ are index-mates, where $\mathcal I$ and $\mathcal J$ are the sets containing the indices of the effective actuators and effective sensors in $\bar B$ and $\bar C$, respectively, and zero otherwise. Now, consider the failure of $k$ communication links corresponding to the $k$ pairs of index-mates $\left (i_{\alpha}, j_{\beta} \right)$, associated with effective dedicated actuators/sensors that control/ the same state variable. Then, by similar reasoning in the proof of Theorem~\ref{theorem_SSSFM}, the only functioning communication link between the pair $\left (i_{\alpha}, j_{\beta} \right)$ will ensure that the tuple $\left( {\bar E, \bar A, \bar B, \bar C; \bar K} \right)$ does not have SSSFM. Therefore, to ensure feasibility, $k+1$ communication channels must be established between the $k+1$ pairs of the index-mates $\left (i_{\alpha}, j_{\beta} \right)$ for each of the pivot. Furthermore, $\bar K$ is also optimal, since each pair of index-mates requires exactly one communication link. Hence, $(\bar B,\bar C,\bar K)$ ensures the robustness of the SSSFM property of the system  $\left( {\bar E, \bar A, \bar B, \bar C; \bar K} \right)$ w.r.t. $k$ failures of \mbox{actuation-sensing-communication} channels.\hfill $\blacksquare$\\
%
%
\textbf{Proof of Theorem~\ref{robustSSSCcost}}: Let the stair matrix $\bar M=P_r\bar A^{\lambda}P_c$ contain $k'$ normalizable step differences $\{\Delta^s\}_{s=1}^{k'}$, where $\bar A^{\lambda}=\bar A-\lambda \bar E$, $P_r$ and $P_c$ are the permutation matrices with appropriate dimensions, and $\Delta^s \in \{0,\times,\otimes \}^{m_1^s\times m_2^s}$. Consider the submatrix $W^B_s \in \mathbb{R}^{m_1^s\times m_2^s}$ of $P_rW^BP_c$, whose indices are the same w.r.t. $\bar M$. For the given collection of step differences  $\{\Delta^s\}_{s=1}^{k'}$, select as a pivot the entry that corresponds to the maximum value $W^B_s\left[p_r^s, p_c^s\right]$ in the submatrix $W^B_s$, where $p_r^s$ and $p_c^s$ denote its row and column indices (w.r.t. $\bar A^\lambda$). Let the row and column indices in $W_s^B$ be $\{p_r^s,r^s_1,\ldots,r^{s}_{m_1^s-1}\}$ and $\{p_c^s,c^s_1,\ldots,c^{s}_{m_2^s-1}\}$, respectively. By selecting an entry $W^B_s\left[r^s_j,\gamma^s_j\right]$ in the row $W^B_s\left[r^s_j,:\right]$ (since the cost of actuating a state does not depend on the actuator used), where $j \in \{1,\dots, m_1^s-1 \}$ and $\gamma^s_j \in \{c^s_1,\ldots,c^{s}_{m_2^s-1}, p_c^s\}$, results in $\bar B_i=[\mathbb{I}^c_n(\mathcal I_i) \ \mathbf{0}_{n\times (\left(k+1\right)n-|\mathcal I_i|)}]$ (up to permutations of columns) which is in the feasibility space of the solutions to $\mathcal P_1^0$. It follows that $\bar B_i$ will incur minimum actuation cost $\sum\limits_{s = 1}^{k'} {\sum\limits_{j = 1}^{m_1^s - 1} {W_s^B} } \left[ {r_j^s,\gamma _j^s} \right]$, where $\mathcal I_i=\bigcup_{s=1}^{k'}\{r^s_1,\ldots,r^{s}_{m_1^s-1}\}$, or $\mathcal I_i=\{1,\ldots, n\}\setminus(\bigcup_{s=1}^{k'}p_r^s)$. As a consequence, consider all possible sequence of collection of the step differences, and select the one with the maximum sum of the pivots' cost described by $W^B$. In addition, robustness can be ensured by following the same reasoning as that presented in the proof of Theorem~\ref{robustSSSC}. \hfill $\blacksquare$\\
%
%
\textbf{Proof of Theorem~\ref{theorem_SSSFMcost}}:  Proof follows similar steps as in Theorem~\ref{robustSSSCcost}. The solution to $\mathcal P_3$ requires selection of pivots which will maximize actuation-sensing-communication cost w.r.t. the sub matrix $W_s$ of $P_rWP_c$ whose indices are the same w.r.t. $\bar M$, where $W=W^B+W^C+W^K$ and $W_s \in \mathbb{R}^{m_1^s\times m_2^s}$. In addition, due to the uniqueness of the cost associated with the index mates, one has to determine the smallest subcollection of these that incur in the minimum cost, which exists among the possible alternatives due to the choice of pivots incurring in the maximum cost. Finally, the robustness can be achieved by invoking the same reasoning as that presented in the proof of Theorem~\ref{theorem_SSSFMrobust}.  \hfill $\blacksquare$ 

{

\appendix

Assumption 1 plays an important role in proving the necessity of Theorem~\ref{SSSCthm}. Specifically, in the latter we can read ``The necessity of the proposed criterion follows similar arguments to those presented in proof of Theorem~1 in~\cite{camsap}.'' It turns out that this argument states the existence of a set of parameters for which the vectors in the step difference are linearly dependent. Nonetheless, this should be 
 ``for \text{all} possible parameters of one vector in the step difference,  the remaining vectors in the step difference admit a parameterization that makes all vectors proportional to each other.'' 

To illustrate how Assumption~1 would lead to the optimal number of dedicated inputs, consider 
\[
\bar A^{\lambda}=\left[\begin{array}{ccc}
\times & \otimes & \otimes\\
\otimes & \otimes & \otimes\\
\otimes &\otimes & \otimes
\end{array}
\right],
\] 
and by invoking Theorem~\ref{dedicatedSol} we  need 
\[
\bar B=\left[\begin{array}{cc}
0 & 0\\
\times & 0\\
0 & \times
\end{array}
\right],
\] 

such that there exist a ramp matrix for $[\bar A^{\lambda} \ \bar B]$ given by
\[
\left[\begin{array}{ccccc}
\otimes & \otimes & \times& 0 &0\\
\otimes & \otimes & \otimes & \times & 0\\
\otimes &\otimes & \otimes &0 & \times
\end{array}
\right],
\] 
since otherwise it is easy to see that having all $\otimes$ set to zero will still ensure the rank to be equal to three.

Now, let us consider one scenario where Assumption~1 does not hold. Suppose we have the following system
\[
\bar A^{\lambda}=\left[\begin{array}{ccc}
\otimes & 0 & \times\\
\times & \otimes & 0\\
\times &\times & \otimes
\end{array}
\right],
\] 
and
\[
\bar B=\left[\begin{array}{cc}
\times\\
0\\
0
\end{array}
\right],
\] 
then the system is selective strong structural controllable (SSSC) since
\[
rank \left(\left[\begin{array}{cccc}
\otimes & 0 & \times&\times\\
\times & \otimes & 0&0\\
\times &\times & \otimes &0
\end{array}
\right]\right)=3
\]
for all possible choices of parameters, but $[\bar A^{\lambda} \ \bar B]$ is not a ramp matrix. If we invoked Theorem~\ref{dedicatedSol}, the conclusion would be  
\[
\bar B=\left[\begin{array}{cc}
0 & 0\\
\times & 0\\
0 & \times
\end{array}
\right],
\] 
ensures  $(\bar A^{\lambda}, \bar B)$ to be SSSC since  $[\bar A^{\lambda} \ \bar B]$ is a ramp matrix, i.e.,
\[
\left[\begin{array}{ccccc}
\otimes & 0 & \times&0&0\\
\times & \otimes & 0&\times&0\\
\times &\times & \otimes &0&\times
\end{array}
\right].
\]

Consequently, if Assumption~1 is not fulfilled, the solutions to the proposed problems are suboptimal (see corrigendum issued for~\cite{popli2019selective}, i.e., the published version of this manuscript). Lastly, it is important to remark that the IEEE 5-bus system explored in this paper satisfies Assumption~1, and therefore, the design and the solution obtained is optimal. Lastly, the problem of obtaining the minimum number of dedicated inputs to ensure SSSC is NP-hard~\cite{Trefois15,Monshizadeh14}, whereas under Assumption~1 it reduces to obtaining a ramp structure which can be done in polynomial-time~\cite{camsap}. Therefore, it would be interesting to explore other assumptions that allow solutions to the proposed problems in polynomial time. 

\section*{Acknowledgement}
We offer our sincere gratitude  to Professor M. Kanat Camlibel (and his co-authors Jiajia Jia, Henk J. van Waarde, and Harry L. Trentelman) for raising relevant question on necessity of ramp structure for ensuring SSSC.

}


\bibliographystyle{IEEEtran}
\bibliography{IEEEabrv,RefTAC}


\end{document}